\documentclass{article}
\pdfoutput=1

\usepackage[T1]{fontenc}  
\usepackage[utf8]{inputenc}

\usepackage[a4paper]{geometry}

\usepackage{amsmath,amsfonts}
\usepackage{stmaryrd} 
\usepackage[pdftex]{graphicx}

\usepackage{authblk}

\usepackage{algorithm}
\usepackage{algpseudocode}

\usepackage{lineno}
\newcommand*\patchAmsMathEnvironmentForLineno[1]{%
  \expandafter\let\csname old#1\expandafter\endcsname\csname #1\endcsname
  \expandafter\let\csname oldend#1\expandafter\endcsname\csname end#1\endcsname
  \renewenvironment{#1}%
     {\linenomath\csname old#1\endcsname}%
     {\csname oldend#1\endcsname\endlinenomath}}%
\newcommand*\patchBothAmsMathEnvironmentsForLineno[1]{%
  \patchAmsMathEnvironmentForLineno{#1}%
  \patchAmsMathEnvironmentForLineno{#1*}}%
\AtBeginDocument{%
\patchBothAmsMathEnvironmentsForLineno{equation}%
\patchBothAmsMathEnvironmentsForLineno{align}%
\patchBothAmsMathEnvironmentsForLineno{flalign}%
\patchBothAmsMathEnvironmentsForLineno{alignat}%
\patchBothAmsMathEnvironmentsForLineno{gather}%
\patchBothAmsMathEnvironmentsForLineno{multline}%
}

\usepackage{booktabs}
\usepackage{color}
\usepackage{siunitx}
\usepackage[numbers,sort]{natbib}
\newcommand{\dd}{\mathrm{d}}
\newcommand{\AT}{\mathit{\psi}}

\newcommand{\tG}{\mathbf{G}}

\newcommand{\vf}{\mathbf{F}}
\newcommand{\tI}{\mathbf{I}}

\newcommand{\sml}{\sigma_\mathrm{l}}
\newcommand{\smt}{\sigma_\mathrm{t}}

\newcommand{\tGbar}{\tG_\mathrm{m}}

\newcommand{\Vm}{V_\mathrm{m}}

\newcommand{\E}{\mathrm{E}}
\newcommand{\Var}{\mathrm{Var}}
\newcommand{\Cov}{\mathrm{Cov}}

\DeclareMathOperator{\dev}{dev}

\newcommand{\bs}{\mathbf}

\title{High-dimensional and higher-order \\ multifidelity Monte Carlo estimators}
\author[1]{A.~Quaglino}
\author[1]{S.~Pezzuto}
\author[1]{R.~Krause}
\affil[1]{Center for Computational Medicine in Cardiology,\par
Institute of Computational Science, \par
Universit\`a della Svizzera italiana,
Lugano, Switzerland}

\date{\small\sffamily Last update: \today}

\begin{document}

\maketitle

\begin{abstract}
Multifidelity Monte Carlo methods rely on a hierarchy of possibly less accurate but statistically correlated
simplified or reduced models, in order to accelerate the estimation of statistics of high-fidelity models without 
compromising the accuracy of the estimates. This approach has recently gained widespread attention in uncertainty quantification \cite{peherstorfer2018}. This is partly due to the availability of optimal strategies for the estimation of the expectation of scalar quantities-of-interest \cite{peherstorfer2015optimal}. In practice, the optimal strategy for the expectation is also used for the estimation of variance and sensitivity indices \cite{qian2017multifidelity}. However, a general strategy is still lacking for vector-valued problems, nonlinearly statistically-dependent models, and estimators for which a closed-form expression of the error is unavailable. The focus of the present work is to generalize the standard multifidelity estimators to the above cases. The proposed generalized estimators lead to an optimization problem that can be solved analytically and whose coefficients can be estimated numerically with few runs of the high- and low-fidelity models. We analyze the performance of the proposed approach on a selected number of experiments, with a particular focus on cardiac electrophysiology, where a hierarchy of physics-based low-fidelity models is readily available.
\end{abstract}

\section{Introduction}
Sampling methods for Uncertainty Quantification (UQ) require multiple evaluations of the problem at hand. 
Often, each sample requires computing the solution of partial differential equations (PDEs), which is generally
a computationally demanding task. For example, in cardiac electrophysiology, a single patient-tailored simulation 
based on the bidomain model can take thousands of node-hours
on a large cluster~\cite{niederer2011}. Therefore, UQ 
for such models at organ scale is currently unfeasible with plain Monte Carlo methods, although its importance has been highlighted in the literature \cite{pathmanathan2017applicability}.
A more sophisticated approach is to perform most of the simulations on a
hierarchy of low-resolution models, yielding the so-called (geometric) multilevel Monte
Carlo (MLMC) method~\cite{gilesmultilevel}, which extends the idea of control
variates methods~\cite{fishman_monte_2003}. The use of the
highest-resolution level guarantees convergence, while a significant portion of
the computational load is offset to the low-resolution hierarchy. While this
can offer a significant speedup of the estimation, it is not always possible to
coarsen complex geometries and to robustly transfer
information from one level to another, to an extent that the full
potential of this approach can be reached~\cite{biehler2015towards}.  This is
for example the case in the context of cardiac
modeling~\cite{rouchdy2007}. Moreover, the estimation of high-order moments requires
special care~\cite{pisaroni2017quantifying,bierig2015convergence}.

An alternative approach, which does not rely on geometry coarsening, is to use model reduction strategies, 
such as projection-based or surrogate models~\cite{gian, corrado2015identification}. 
Such an approximation is build upon the observation that in many applications it is clearly possible to
distinguish between \emph{offline} and \emph{online} phases of the workflow,
where the former typically can be very expensive and encompass everything that
can be precomputed in advance, e.g., training a surrogate model before any
patient-specific data becomes available, while the latter should be very cheap and
consist only of evaluations of the reduced model. Unfortunately, this approach has two
significant drawbacks. On the one hand, complex error estimates must be
provided to ensure that the approximation error is within acceptable bounds. On
the other hand, repeating the offline training is very expensive but may often be
necessary, e.g., in the case of patient-specific UQ. Alternative to that, physics-based 
reduction strategies are also possible: in cardiac electrophysiology, the activation map
can readily be computed via the eikonal equation, which provides a physiologically-meaningful solution~\cite{pullan2002}, 
strongly linked to the bidomain equation~\cite{colli93}. Error bounds, however,
are not trivially found, limiting the applicability of the methodology.

A more recent idea is to replace the above paradigm of model \emph{selection} with that of model \emph{fusion}, 
where low-fidelity models are allowed to be inaccurate, in the sense of not providing
an approximation within certain error bounds, provided that they exhibit some degree of
correlation~\cite{peherstorfer2015optimal} or even just statistical
dependence~\cite{koutsourelakis2009accurate} to the high-fidelity one. This
approach has a twofold advantage. Firstly, it is the statistical
dependence, rather than the error bounds of coarse models, that is crucial to
ensure that propagating uncertainties via the low-fidelity models provides
useful information on the statistics of the high-fidelity quantity-of-interest.
Secondly, in the case of PDEs, low-fidelity models are not restricted to low-resolution
geometries, therefore solving them can be several orders of magnitude faster
than a coarse model. If a statistical dependence between the models is
present, in practice, most of the computational effort 
is taken by the estimation of the linear correlation coefficients or a nonlinear map between the outputs of the low- and the high-fidelity models.
For such reasons, multifidelity Monte Carlo (MFMC) methods
have become very popular over the last years and their applications span the
fields of UQ, inverse problems, and optimization~\cite{peherstorfer2018}.

Several approaches can be used to combine the outputs of the model hierarchy. The use of a Bayesian regression,
mapping the low-fidelity output to the high-fidelity one, has been advocated in~\cite{koutsourelakis2009accurate}. In \cite{Quaglino2018}, a Gaussian Process Regression is considered in the context of cardiac electrophysiology. The Bayesian nature of this
approach automatically augments the estimate with full probability distributions and credible intervals,
obtained independently of the degree of statistical dependence between models and the number of samples
employed. This makes it a very suitable approach for scenarios where resources are scarce, such as 
clinical applications. In contrast to that, the multifidelity Monte Carlo method~\cite{peherstorfer2015optimal} 
provides point estimates, rather than probability distributions. However, these are equipped with convergence estimates 
for the asymptotic case, i.e., their error vanishes as the number of samples goes to infinity. This makes this approach 
particularly interesting when computational resources are large (but finite) and accurate
estimates are desired.

Current MFMC methods have been successfully used for estimating the expectation,
the variance, the Sobol indices, and rare events~\cite{peherstorfer2018}. For the expectation
of a scalar quantity-of-interest, the optimal multifidelity estimator for a given computational budget can be found analytically~\cite{peherstorfer2015optimal}.
Similarly, \citet{qian2017multifidelity} proposed an estimator for variance
and Sobol indices, with an analytical derivation only for the former.
Interestingly,
the authors also report that the optimal weights and sampling strategy for the
expectation multifidelity estimator
performs very well also for the variance and Sobol indices, although not being optimal.
This observation suggests that a multifidelity estimator for a general statistic
of interest may perform relatively well even when the weights and sampling strategy
are not optimal.

In this work, we explore the above idea, hence extending the aforementioned multifidelity estimators in several ways.
First, we generalize the scalar-valued multifidelity estimator to vector-valued quantity-of-interests,
proposing an optimal estimator in terms of mean squared error (MSE) that closely resembles the one
obtained in the scalar case. Second, we analyze a broad class of high-order statistics
of interest, for which a closed-form expression of the error is unavailable and therefore the optimal sampling strategy is unknown. To overcome this limitation, we suggest to approximate the MSE by estimating the correlations between the estimators, via sampling of the low- and high-fidelity models. Third, we suggest a strategy to leverage nonlinear statistical dependence between the models, which is closely connected to \cite{koutsourelakis2009accurate,Quaglino2018}.

The paper is organized as follows: Section 2 introduces the standard multifidelity Monte Carlo estimators of scalar quantities-of-interest. Section 3 describes the proposed generalized multifidelity estimators. Section 4 presents the results of the numerical experiments for the assessment of each estimator. Section 5 details a real-world application to cardiac electrophysiology. Section 6 illustrates an application to nonlinear incompressible elasticity.

\section{Standard multifidelity Monte Carlo estimators}
\label{sec:methods}

Let upper-case letters denote random
variables (e.g.\ $\bs{S}$, $\Psi$) and lower-case letters the values these take
(e.g.\ $\bs{s}$, $\psi$), where $\psi$ and $\psi(\bs{s})$
denote both, the value of the quantity-of-interest (QoI) as well as the function that provides the
QoI with respect to the random input $\bs{s}$.
Let $\psi^{(i)}(\bs{s})$ be the $i$-th model of the multifidelity hierarchy, with $i=1,\ldots,K$. For ease of notations, we
denote by $\psi_h(\bs{s}) := \psi^{(1)}(\bs{s})$ and $\psi_l(\bs{s}) := \psi^{(K)}(\bs{s})$ the high- and lowest-fidelity models.

The main task in the context of uncertainty propagation is the computation
of integrals. In particular, the sought expectation of the QoI is obtained by:
\begin{equation}
\E[\Psi_h] = \int_{\Omega} \psi_h(\bs{s}) ~p(\bs{s})\:\dd\bs{s},
\label{eq:expect}
\end{equation}
where $p(\bs{s})$ denotes the probability density of the random inputs.
Similarly, one can compute other statistical indicators
(e.g. moments, densities etc), by evaluating integrals of the form $\int
h(\psi_h(\bs{s})) ~p(\bs{s})ds$. The method of choice in high-dimensional settings is
direct Monte Carlo, where the expectation of Equation \eqref{eq:expect} is
estimated by:
\begin{equation}
	\hat{\Psi}^{(i)}_m =\frac{1}{m} \sum_{j=1}^m \psi^{(i)}(\bs{s}_j),
	\label{eq:mc}
\end{equation}
where $i=1$ and the $\bs{s}_j$ are independent, identically distributed samples drawn from $p(\bs{s})$. 
The convergence rate is independent of the dimension of 
$\bs{S}$ and the error decays as $\mathcal{O}(\frac{1}{\sqrt{m}})$ 
\cite{kalos_monte_2008}.
Therefore, in problems where each evaluation of $\psi_h(\bs{s})$ poses a 
significant computational burden, the use of direct Monte Carlo can become 
impractical or even infeasible.

\subsection{Scalar multifidelity expectation estimator}
At this stage, a multifidelity estimator is introduced
\begin{equation}
	\hat{\Psi}_h = \hat{\Psi}^{(1)}_{m_1} + \sum_{i=2}^K \alpha_i \left( \hat{\Psi}^{(i)}_{m_i} - \hat{\Psi}^{(i)}_{m_{i-1}} \right).
	\label{eq:mfmc}
\end{equation}
The key aspect of the multifidelity estimate is that the magnitude of the discrepancy
$|\psi_h(\bs{s}) -  \psi^{(i)}(\bs{s})|$ does not affect its accuracy. More precisely, let 
the Pearson correlation coefficient of the $i$-th model be
\begin{equation}
	\rho_{i,j} := \frac{\Cov[\psi^{(i)},\psi^{(j)}]}{\sigma_i \sigma_j}
	\label{eq:pearson}
\end{equation}
and its variance $\sigma^2_i := \Var[\hat{\Psi}^{(i)}]$. Then, it is possible to show \cite{peherstorfer2015optimal} that the MSE of the multifidelity estimator is
\begin{equation}
	\Var[\hat{\Psi}_h] = \frac{\sigma^2_1}{m_1} + \sum^K_{i=2} \left( \frac{1}{m_{i-1}} - \frac{1}{m_i} \right) \left( \alpha_i^2 \sigma^2_i -2 \alpha_i \rho_{1,i} \sigma_1 \sigma_i \right).
	\label{eq:variance}
\end{equation}
Moreover, the cost of the estimator is 
\begin{equation}
	C(\hat{\Psi}_h) = \sum^K_{i=1} w_i m_i,
\end{equation}
where $w_i$ is the computational cost of the $i$-th model. Under some mild assumptions on $\rho_{i,j}$ and $w_i$, it is possible to perform a constrained minimization of the variance \eqref{eq:variance} with respect to $\alpha$ and $m$, for a prescribed computational cost $C(\hat{\Psi}_h) = B > 0$, obtaining
\begin{equation}
	\alpha^*_{i} := \frac{\rho_{1,i}\sigma_1}{\sigma_i}, \quad m^*_1=\frac{B}{\sum_{i=1}^K w_i r_i}, \quad m^*_i=m^*_1 r_i, \quad i=2,\ldots,K
	\label{eq:alpha}
\end{equation}
where 
\begin{equation}
	r_{i} := \sqrt{\frac{w_1\left(\rho^2_{1,i}-\rho^2_{1,i+1}\right)}{w_i\left(1-\rho^2_{1,2}\right)}}.
	\label{eq:t}
\end{equation}
It follows that the ratio of errors (i.e. the variance reduction) of the standard and the multifidelity estimators, for the same computational budget, is 
\begin{equation}
	\frac{E(\hat{\Psi}^*_{h})}{E(\hat{\Psi}^{(1)}_n)} = \left( \sum_{i=1}^K \sqrt{\frac{w_i}{w_1}\left(\rho^2_{1,i}-\rho^2_{1,i+1}\right)}\right)^2,
	\label{eq:ratio}
\end{equation}
where $n=B/w_1$. Therefore, the variance is reduced if the costs $w_i$ and the differences of the squared correlation coefficients are low.
This condition replaces the classical assumption on the deterministic pointwise errors $|\psi_h(\bs{s}) -  \psi^{(i)}(\bs{s})|$ of the low-fidelity models with respect to the high-fidelity one.

\subsection{Variance and Sobol indices estimators}
The standard approach to develop multifidelity estimators for higher-order moments $q(\Psi_h)$ is to substitute $\hat{\Psi}^{(i)}_m$ in \eqref{eq:mfmc} with the unbiased single-level estimator $\hat{q}$. For example, in the case of variance and Sobol indices, the use of the following estimators is proposed in \cite{qian2017multifidelity}
\begin{align*}
	\hat{V} &= \frac{1}{m-1} \sum_{i=1}^m \left( \psi(s_i) - \hat{\Psi} \right)^2, \\
	\hat{V_j} &= \frac{2}{2m-1} \left( \sum_{i=1}^m \psi(s_i) \psi(y^j_i) - m \left( \frac{\hat{\Psi}+ \hat{\Psi}'}{2} \right)^2 + \frac{\hat{V}+ \hat{V}'}{4} \right), \\
	\hat{T_j} &= \frac{1}{2m} \sum_{i=1}^m \left( \psi(s'_i) - \psi(y^j_i) \right)^2,
\end{align*}
where $s'$ is a second set of samples and $y^j_i$ is equal to $s'_i$, except for the $j$-th component which is equal to that of $s_i$.  However, obtaining error estimates in the form of \eqref{eq:variance} is generally too difficult or impractical. Therefore, it is advocated in \cite{qian2017multifidelity} to employ the optimal strategy $(m^*,\alpha^*)$ for the expectation also for the computation of other estimators. In the following section, we propose a different approach based on estimating the correlations between each pairs of summands (as computed via two different models) in the above estimators.

\section{Generalized multifidelity Monte Carlo estimators}

\subsection{Vector-valued multifidelity expectation estimator} \label{vectorSec}
We now consider the case of a vector-valued quantity-of-interest. In particular, we focus on a spatially-dependent $\psi$, i.e. $\psi(x)$. It follows that the random variables and their estimators are also spatially dependent, i.e., $\Psi(x)$ and $\hat{\Psi}(x)$. We consider the finite-dimensional case of computing them at a given set of points $\{ x_j \}_{j=1}^N$. Given a computational budget $B$, the optimal scalar-valued estimator can therefore be computed at each point $x_j$. By prescribing a tolerance $\epsilon^2$ for the integrated MSE of the optimal estimator, expressed as a function of the computational budget $B$, as given in \cite{peherstorfer2015optimal}, we obtain
\begin{equation*}
	\int_{\Omega} \frac{\sigma_1^2(x)}{B} \left( \sum_{i=1}^K \sqrt{\frac{w_i}{w_1}\left(\rho^2_{1,i}(x)-\rho^2_{1,i+1}(x)\right)}\right)^2 dx < \epsilon^2.
\end{equation*}
By assuming a spatially-uniform budget $B(x) = \overline{B}$, linearity of the integral can be exploited, so that the optimal budget $\overline{B}^*$ that satisfies the above equation is easily found. However, this will yield a spatially-dependent number of model evaluations $m_i(x)$. Since it is not possible to perform a model evaluation only for a subset of the points, the following choices are made
\begin{equation*}
	m_i = \max_{j=1, \ldots, N} m_i (x_j).
\end{equation*}
Unfortunately, in general this will lead to a large number of model evaluations. Instead, we define the integrated vector error as 
\begin{align} \label{vectorVal}
	\int_\Omega E(\hat{\Psi}_h,x,\alpha)~ dx := \sum_{j=1}^N E(\hat{\Psi}_h,x_j,\alpha) | \Omega_j |,
\end{align}
where $| \Omega_j |$ is the measure of the local integration area and 
\begin{align*}
	E(\hat{\Psi}_h,x_j,\alpha) &= \frac{\sigma^2_1(x_j)}{m_1} + \sum^K_{i=2} \left( \frac{1}{m_{i-1}} - \frac{1}{m_i} \right) \left( \alpha^2_i(x_j) \sigma^2_i(x_j) -2 \alpha_i(x_j) \rho_{1,i}(x_j) \sigma_1(x_j) \sigma_i(x_j) \right).
\end{align*}
Therefore, the error of the vector-valued case can be seen as a weighted
average of the scalar-valued errors. Hence, the resulting estimator can be also
used for the generic spa\-tially-in\-de\-pendent vector-valued case, given some
importance weights $| \Omega_j |$. Setting the partial derivative of
$E(\hat{\Psi}_h,x_j,\alpha)$ with respect to $\alpha_i(x_j)$ equal to zero, we
obtain that
\begin{align*}
	\alpha^*_{i}(x_j) := \frac{\rho_{1,i}(x_j)\sigma_1(x_j)}{\sigma_i(x_j)}, \quad \forall j=1\cdots N.
\end{align*}
It follows that
\begin{align*}
	E(\hat{\Psi}_h,x_j,\alpha^*) &= \frac{\sigma^2_1(x_j)}{m_1} - \sum^K_{i=2} \left( \frac{1}{m_{i-1}} - \frac{1}{m_i} \right) \rho^2_{1,i}(x_j) \sigma^2_1(x_j).
\end{align*}
Therefore, the error sum can be written as
\begin{align}
	\label{eq:VarianceOptAlpha}
	\sum_{j=1}^N E(\hat{\Psi}_h,x_j,\alpha) = \frac{1}{m_1} \overline{\sigma}^2_1 - \sum^K_{i=2} \left( \frac{1}{m_{i-1}} - \frac{1}{m_i} \right) \ \overline{\rho}_{1,i}^2 \overline{\sigma}_1^2 ,
\end{align}
where
\begin{align*}
	\overline{\sigma}^2_1 := \sum_{j=1}^N \sigma^2_1(x_j) | \Omega_j |, \quad \overline{\rho}^2_{1,i} := \frac{1}{\overline{\sigma}^2_1} \left( \sum_{j=1}^N \rho^2_{1,i}(x_j) \sigma^2_1(x_j) | \Omega_j | \right)
\end{align*}
Equation \eqref{eq:VarianceOptAlpha} has the same form as for the scalar case, hence it leads to
\begin{equation}
	m^*_1=\frac{B}{\sum_{i=1}^K w_i \overline{r}_i}, \quad m^*_i=m^*_1 \overline{r}_i, \quad \overline{r}_{i} := \sqrt{\frac{w_1\left( \overline{\rho}^2_{1,i}-\overline{\rho}^2_{1,i+1}\right)}{w_i\left(1-\overline{\rho}^2_{1,2}\right)}}.
\end{equation}

\subsection{Higher-order and sensitivity estimators}
In general, given a single-level estimator $\hat{q}$, its corresponding multifidelity estimator is defined as
\begin{equation}
	\hat{q}_h = \hat{q}^{(1)}_{m_1} + \sum_{i=2}^K \alpha_i \left( \hat{q}^{(i)}_{m_i} - \hat{q}^{(i)}_{m_{i-1}} \right).
	\label{eq:mfmc_q}
\end{equation}
Regardless of the true expression of $\hat{q}$, we implicitly define a function $q$ via the relation
\begin{equation}
	\hat{q}^{(j)}_{m_i} = \frac{1}{m_i} \sum_{n=1}^{m_i} q(\Psi^{(j)}(s_n)).
	\label{eq:q}
\end{equation}
Under the assumption that the $s_n$ are i.i.d.\ and $m_{i-1} < m_i$, it can be shown that 
\begin{align*}
	\Cov[\hat{q}^{(j)}_{m_i},\hat{q}^{(k)}_{m_{i-1}}-\hat{q}^{(k)}_{m_i}]= \frac{1}{m_i}  \left( 1 - \min \left( 1, \frac{m_i}{m_{i-1}} \right) \right) c_{jk} =0, \quad \forall (j,k), 
\end{align*}
where $c_{jk} = \Cov[q(\Psi^{(j)}),q(\Psi^{(k)})]$. Therefore, it possible to write the MSE of \eqref{eq:mfmc_q} as
\begin{align*}
	\Var[\hat{q}_h] = \Var[\hat{q}^{(1)}_{m_1}] + &\sum_{i=2}^K \alpha_i^2 \left( \Var[\hat{q}^{(i)}_{m_{i-1}}] - \Var[\hat{q}^{(i)}_{m_{i}}]  \right) \\
	- 2 &\sum_{i=2}^K \alpha_i \left( \Cov[\hat{q}^{(1)}_{m_1},\hat{q}^{(i)}_{m_{i-1}}] - \Cov[\hat{q}^{(1)}_{m_1},\hat{q}^{(i)}_{m_i}] \right).
\end{align*}
While in general it is not possible to find analytic expressions for the above terms, they can be estimated numerically in a preprocessing step by the definition of $q$ in equation \eqref{eq:q} as 
\begin{align*}
	\Cov[\hat{q}^{(j)}_{m_i},\hat{q}^{(k)}_{m_{i-1}}] \approx \frac{1}{m_i} \frac{1}{N-1} \sum_{n=1}^N \left( q(\Psi^{(j)}(s_n)) - \hat{q}^{(j)}_N \right) \left(q(\Psi^{(k)}(s_n)) - \hat{q}^{(k)}_N \right), \quad m_{i-1} \leq m_i, 
\end{align*}
where $N$ is a small number (i.e. $N \approx 10-20$). By defining 
\begin{align*}
	\left( \sigma^q_i \right)^2 & := \frac{1}{N-1} \sum_{n=1}^N \left( q(\Psi^{(i)}(s_n)) - \hat{q}^{(i)}_N \right)^2, \\
	\rho^q_{1,i} & := \frac{1}{\sigma^q_1 \sigma^q_i (N-1)} \sum_{n=1}^N \left( q(\Psi^{(1)}(s_n)) - \hat{q}^{(1)}_N \right) \left(q(\Psi^{(i)}(s_n)) - \hat{q}^{(i)}_N \right),
\end{align*}
we obtain an expression similar to \eqref{eq:variance}
\begin{equation}\label{eq:mfmc_Jq}
	\Var[\hat{q}_h] = \frac{(\sigma^q_1)^2}{m_1} + \sum^K_{i=2} \left( \frac{1}{m_{i-1}} - \frac{1}{m_i} \right) \left( \alpha_i^2 (\sigma^q_i)^2 -2 \alpha_i \rho^q_{1,i} \sigma^q_1 \sigma^q_i \right).
\end{equation}
Therefore, the optimization procedure presented above leads to a similar result with $\sigma^q$ and $\rho^q$ in place of $\sigma$ and $\rho$.

\subsection{Nonlinearly statistically-dependent estimators}
We now consider the case where the correlation coefficient $\rho_{1,i}$ between the models is low, but the two model exhibit a high degree of nonlinear statistical dependence. For ease of notation, we consider to estimate the expectation of a scalar $\Psi$, but the same procedure can be applied to estimate any statistics and to the vector-valued case, by combining it with the strategies from the previous sections. Let us consider an estimator of the form 
\begin{equation}
	\hat{g}_h = \hat{\Psi}^{(1)}_{m_1} + \sum_{i=2}^K \alpha_i \left( \hat{g}^{(i)}_{m_i} - \hat{g}^{(i)}_{m_{i-1}} \right),
	\label{eq:mfmc_g}
\end{equation}
where 
\begin{equation}
	\hat{g}^{(j)}_{m_i} = \frac{1}{m_i} \sum_{n=1}^{m_i} g(\Psi^{(j)}(s_n)).
	\label{eq:g}
\end{equation}
The function $g$ is an additional 1-dimensional model, for example a Gaussian Process regression, which is fitted using $N_g$ samples. The idea is that for a strongly nonlinear statistical dependence, $N_g+N_{\rho} \ll N$, where $N$ and $N_{\rho}$ are the number of samples used to estimate the correlation for, respectively, the estimators \eqref{eq:mfmc} and \eqref{eq:mfmc_g}. Following the same procedure as in the previous section, it follows that
\begin{equation}\label{eq:mfmc_Jg}
	\Var[\hat{g}_h] = \frac{(\sigma_1)^2}{m_1} + \sum^K_{i=2} \left( \frac{1}{m_{i-1}} - \frac{1}{m_i} \right) \left( \alpha_i^2 (\sigma^g_i)^2 -2 \alpha_i \rho^g_{1,i} \sigma_1 \sigma^g_i \right).
\end{equation}
Therefore, the optimization procedure presented above leads to a similar result with $\sigma^q$ and $\rho^q$ in place of $\sigma$ and $\rho$, where
\begin{align*}
	\left( \sigma^g_i \right)^2 & := \frac{1}{N-1} \sum_{n=1}^N \left( g(\Psi^{(i)}(s_n)) - \hat{g}^{(i)}_N \right)^2, \\
	\rho^g_{1,i} & := \frac{1}{\sigma_1 \sigma^g_i (N-1)} \sum_{n=1}^N \left( \Psi^{(1)}(s_n) - \hat{\Psi}^{(1)}_N \right) \left(g(\Psi^{(i)}(s_n)) - \hat{g}^{(i)}_N \right).
\end{align*}
%

\subsection{Choosing the computational budget}
In order to choose the budget, we express the integrated MSE of the optimal estimator as a function of the computational budget $B$, as given in \cite{peherstorfer2015optimal}
\begin{equation*}
	E(\hat{\Psi}^*_{h},B) = \frac{\overline{\sigma}_1^2}{B} \left( \sum_{i=1}^K \sqrt{\frac{w_i}{w_1}\left(\overline{\rho}^2_{1,i}-\overline{\rho}^2_{1,i+1}\right)}\right)^2 < \epsilon^2,
\end{equation*}
from which it follows that 
\begin{equation}
	\overline{B}^* = \left( \frac{\overline{\sigma}_1}{\epsilon} \sum_{i=1}^K \sqrt{\frac{w_i}{w_1}\left(\overline{\rho}^2_{1,i}-\overline{\rho}^2_{1,i+1} \right)}\right)^2.
\end{equation}

\section{Numerical assessment}
In this section, we test the proposed generalized estimators on two simple benchmark problems. The first test is a well-known problem in the context of multifidelity methods \cite{qian2017multifidelity}. There, we consider the expectation, variance, and Sobol index estimators. In the second test, we focus on comparing the linear and nonlinearly statistically-dependent estimators. In Figures \ref{ishigami_MV} and \ref{quintic_MV}, the two approaches are referred to as Linear and GPR, respectively. All of the multifidelity routines are implemented in a in-house, newly developed Python library called SLOTH. 

\begin{figure}[!htb]
	\centering
	\includegraphics[width=.9\textwidth]{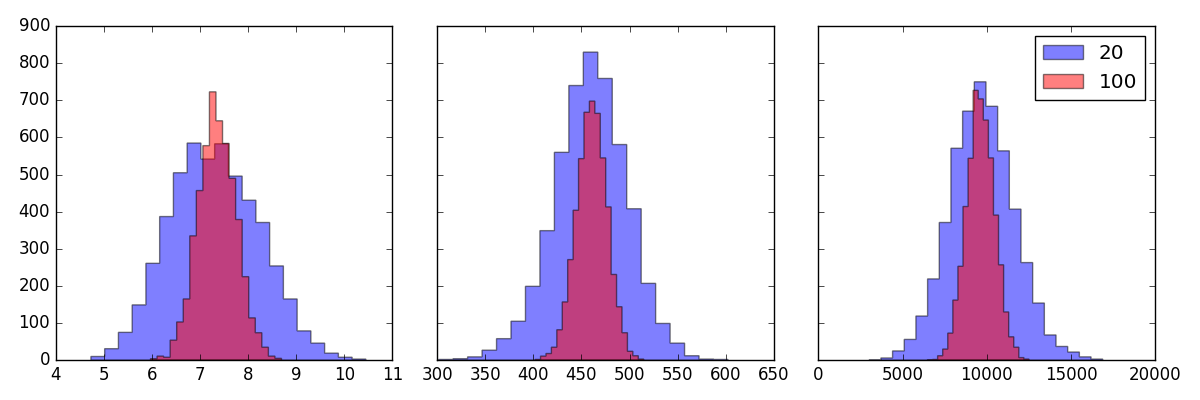}
	\includegraphics[width=.9\textwidth]{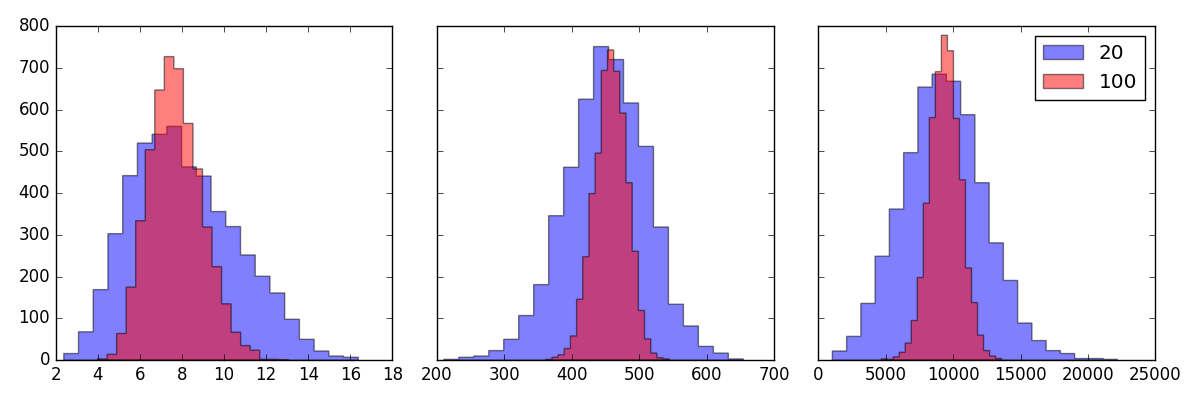}
  \caption{Distribution of optimal $m_k$ for the expectation (top row) and variance (bottom row) estimators for the Ishigami model with 20 (blue) and 100 (red) preprocessing samples, computed with 5000 replicates. Columns: $m_1$, $m_2$, and $m_3$. Plot axes: values (horizontal) and frequency (vertical) of the estimator.}
	\label{ishigami_m_MV}
\end{figure}

\begin{figure}[!htb]
	\centering
	\includegraphics[width=.9\textwidth]{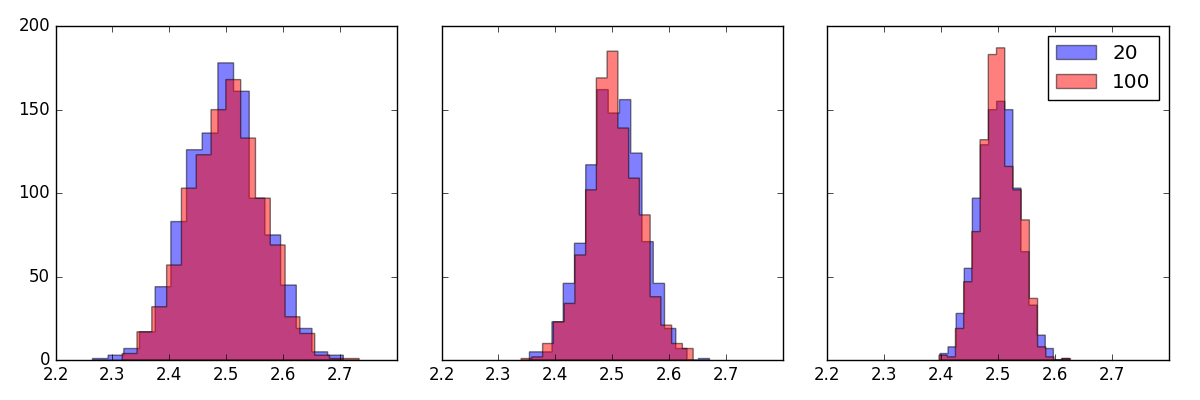}
	\includegraphics[width=.9\textwidth]{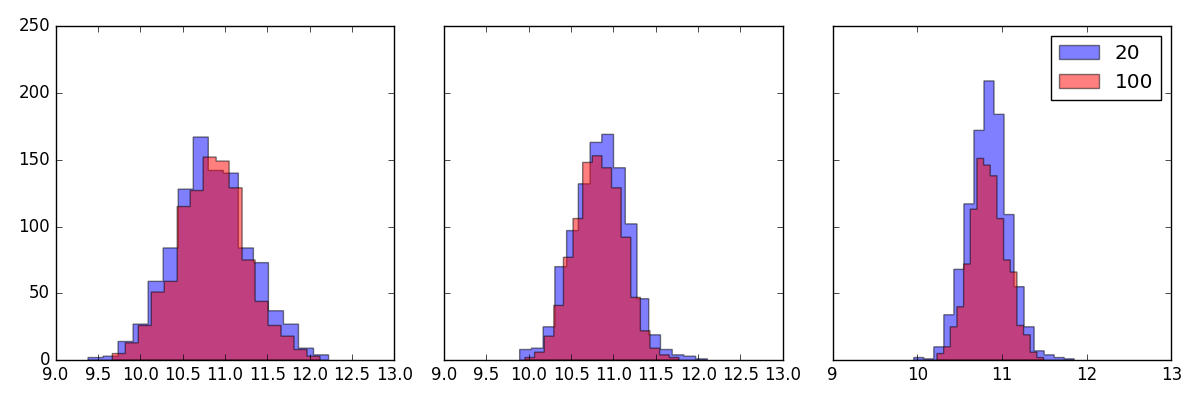}
  \caption{Distribution of the computed estimator of the expectation (top row) and variance (bottom row) for the Ishigami model with 20 (blue) and 100 (red) preprocessing samples, computed with 1000 replicates. Columns: $p=40$, $p=80$, and $p=160$. Plot axes: values (horizontal) and frequency (vertical) of the estimator.}
	\label{ishigami_y_MV}
\end{figure}

\subsection{Ishigami model}
We consider the Ishigami model \cite{ishigami1990importance} and the low-fidelity models proposed in \cite{qian2017multifidelity}
\begin{align*}
	f^{(1)} &=  \sin(z_1) + 5 \sin^2(z_2) +  \frac{1}{10} z_3^4 \sin(z_1), \quad z_i \sim \mathcal{U}(-\pi,\pi) \\
	f^{(2)} &=  \sin(z_1) + 4.75 \sin^2(z_2) + \frac{1}{10} z_3^4 \sin(z_1), \\
	f^{(3)} &=  \sin(z_1) + 3 \sin^2(z_2) + \frac{9}{10} z_3^2 \sin(z_1). 
\end{align*}
For this model, expectation and variance can be computed analytically as $\E[f(z)]=2.5$ and $\Var[f(z)] \approx 10.845$. The reference values of the Sobol indices can be found in \cite{qian2017multifidelity}. We study numerically the errors in the proposed multifidelity estimators. For the expectation and variance estimators, we use both approaches of linear correlation and nonlinear statistical dependence. Moreover, we consider the Sobol indices as a vector-valued output, using the multidimensional estimator proposed in Section~\ref{vectorSec}.

\subsubsection{Preprocessing budget}
Figure \ref{ishigami_m_MV} shows the histograms of the estimated optimal $m$ with 5000 replicates of the preprocessing step. In each replicate, $m$ is obtained by analytic minimization \eqref{eq:alpha} of the error \eqref{eq:mfmc_Jq}, computed with either 20 (blue) or 100
(red) samples. Figure \ref{ishigami_y_MV} shows the histogram of the computed
estimators for 2000 replicates of the full analysis, again performed with
either 20 (blue) or 100 (red) preprocessing samples. From
Figure~\ref{ishigami_m_MV},
it is clear that the number of samples used for preprocessing has an
impact on the accuracy of the computed $m$ (and $\alpha$). This is due to
the accuracy of the variance and correlation estimates, which appear in the
formula for computing $m$ and $\alpha$.
Fortunately, the
second figure shows that this does not have an impact on the computed
estimators. In all of the examples, we have used 100 samples as preprocessing
budget.

\subsubsection{Expectation and variance estimation}
\begin{figure}[!htb]
	\centering
	\includegraphics[width=.9\textwidth]{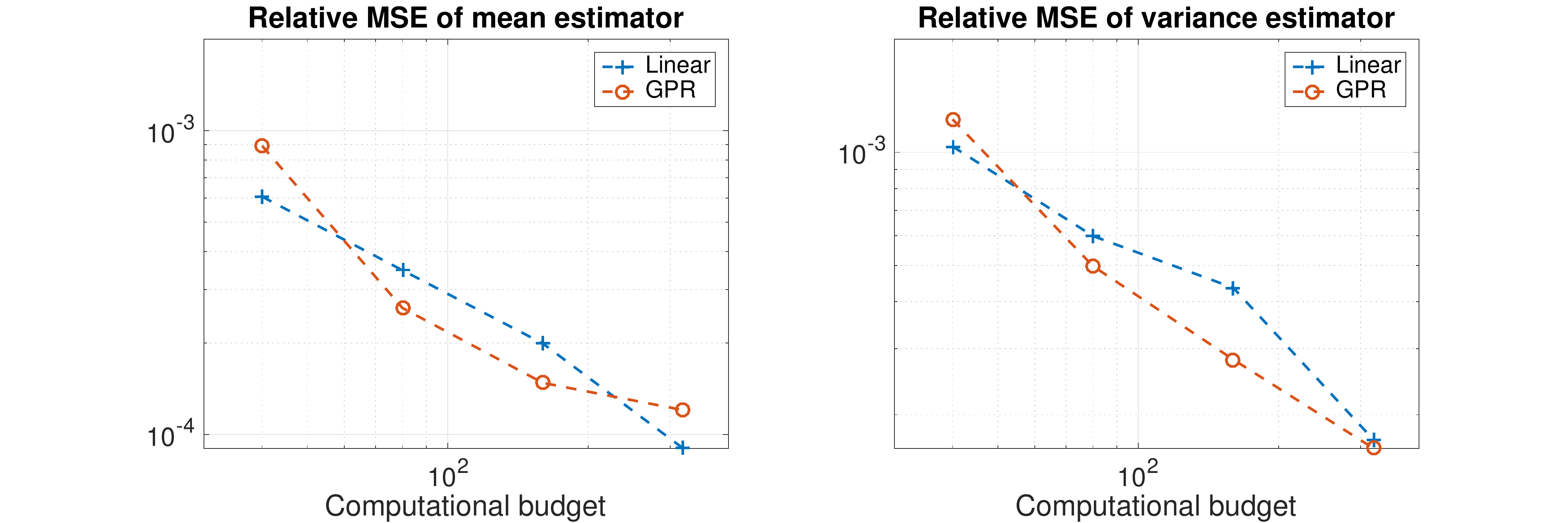}
  \caption{Relative MSE of the proposed expectation and variance estimators for the Ishigami model as a function of computational budget $p$, computed with 100 replicates, for the linear and nonlinear estimators. Blue line: estimator using equation \eqref{eq:mfmc_q}. Red line: estimator using equation \eqref{eq:mfmc_g} and a GPR model is trained with 100 samples.  }
	\label{ishigami_MV}
\end{figure}

Figure \ref{ishigami_MV} shows the computed error of the expectation and variance estimators using 100 replicates of the proposed estimators \eqref{eq:mfmc_q} and \eqref{eq:mfmc_g}. We see that in this case, taking advantage of the nonlinear statistical dependence has no advantage over using only linear correlations. Table \ref{table_ishigami_MV} summarizes the optimal $\alpha$ and $m$ for the case of linear correlations, both averaged over 100 runs. As discussed in \cite{qian2017multifidelity}, the optimal budget allocations for expectation and variance are similar. Therefore, in principle there is no advantage in using a separate allocation for the variance. However, for this problem, we get a better approximation of the optimal $\alpha$ (0.8826 vs. 0.9455, as given in \cite{qian2017multifidelity}). This improves the accuracy of the computed variance with no extra cost, since the samples used for estimating $\rho$ and $\sigma$ in the standard multifidelity estimators, which are required to estimate the expectation, can be re-used for the proposed generalized ones. Moreover, this is achieved without using any a priori error estimate and can therefore be used also in the cases where such bounds are not available.

\begin{table}[!htp]
\centering
\begin{tabular}{lrrrr}
     \toprule
     \multicolumn{1}{l}{Models} &
     \multicolumn{2}{c}{Expectation est.} &
     \multicolumn{2}{c}{Variance est.} \\
     \cmidrule(lr){2-3}
     \cmidrule(lr){4-5}
	 & $m_k$ & $\alpha_k$ & $m_k$ & $\alpha_k$ \\
	 \midrule
	 $f^{(1)}$ & 7  & 1 & 8 & 1 \\
	 $f^{(2)}$ & 461 & 1.0144 & 458 &  1.0144  \\
	 $f^{(3)}$ &  9633 & 0.8826 & 9564 & 0.9289 \\
	 \bottomrule
\end{tabular}
\caption{Estimates for the optimal $\alpha$ and $m$ (with a budget $p=40$) of the Ishigami model obtained via analytic minimization of \eqref{eq:mfmc_Jq} for the case of linear correlations. First two columns: values for the expectation estimator. Last two columns: values for the variance estimator. Both are averaged over 100 runs and estimated using 100 preprocessing samples. }
\label{table_ishigami_MV}
\end{table}

\subsubsection{Sobol indices estimation}
\begin{figure}[!htb]
	\centering
	\includegraphics[width=.9\textwidth]{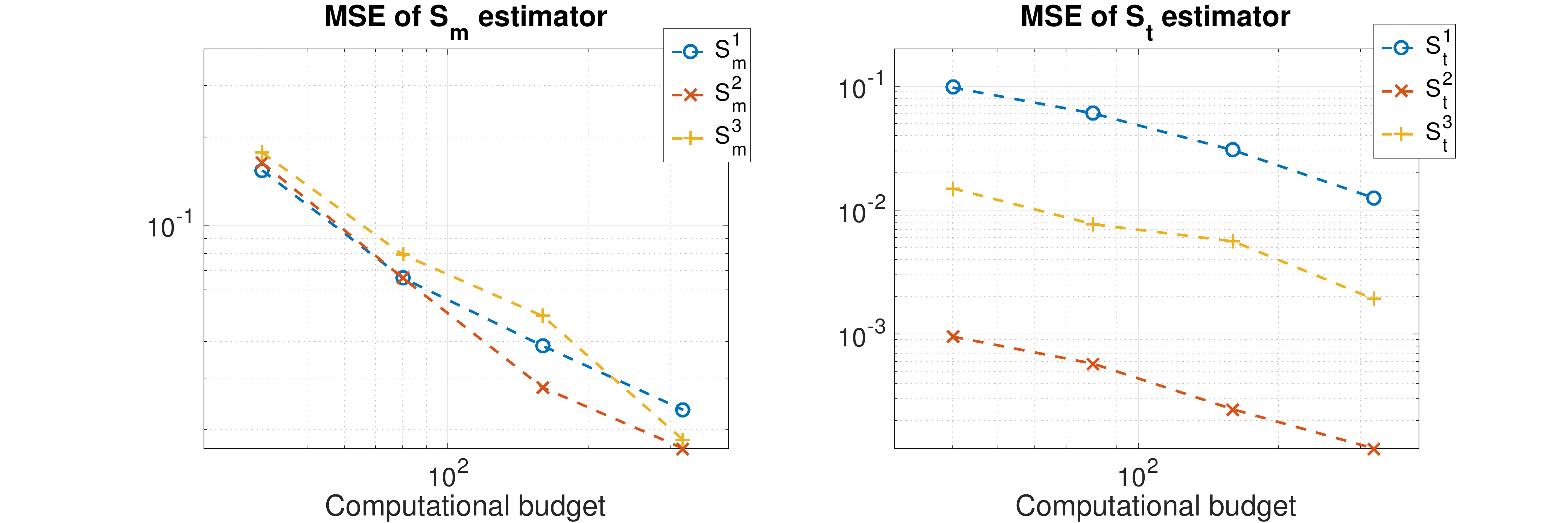}
  \caption{Absolute MSE of the proposed Sobol estimators for the Ishigami model as a function of the computational budget $p$, computed with 100 replicates. Left: main effect sensitivities. Right: total effect sensitivities. Both are computed with the linear vector-valued estimator.}
	\label{ishigami_S}
\end{figure}

Here, we consider the main and total effect Sobol indices as vector-valued outputs of three components each. Therefore, we obtained a single optimal value $m_k$ for each set of sensitivities, as discussed in Section \ref{vectorSec}, rather than a separate optimum for each index. Table \ref{table_ishigami_S} summarizes the estimates for the optimal $\alpha$ and $m$, both averaged over 100 runs. We see that the optimal allocation differs between the two sets, and also with respect to the expectation and variance, shown in Table \ref{table_ishigami_MV}. Therefore, in this case the optimum for the expectation does not guarantee a nearly optimal allocation of the samples. Figure \ref{ishigami_S} shows the computed error of the proposed multifidelity Sobol indices estimators using 100 runs of the multifidelity algorithm, using only linear correlations. In particular, the errors are similar for the main effect sensitivities, while they differ for the total effect ones. We have used uniform weights $| \Omega_j |=1$ in \eqref{vectorVal}, however a different choice would have been possible if a higher precision had been required for one of the Sobol indices.

\begin{table}[htp]
\begin{center}
\begin{tabular}{lrrrrrrrr}
	\toprule
     	\multicolumn{1}{l}{Models} &
     	\multicolumn{4}{c}{Main index est.} &
     	\multicolumn{4}{c}{Total index est.} \\
     	\cmidrule(lr){2-5}
     	\cmidrule(lr){6-9}
	 & $m^M_k$ & $\alpha^M_k$ & $\alpha^M_k$ & $\alpha^M_k$ & $m^T_k$ & $\alpha^T_k$ & $\alpha^T_k$ & $\alpha^T_k$ \\
	 \midrule
	 $f^{(1)}$ & 9  & 1 & 1 & 1 & 0 & - & - & -  \\
	 $f^{(2)}$ & 485 & 1.049 & 1.057 &  1.028  & 55 & 1 & 1.108 & 1 \\
	 $f^{(3)}$ &  7071 & 1.005 & 0.976 & 0.922 & 12471 & 0.828 & 2.778 & 1.051 \\
	 \bottomrule
\end{tabular}
\end{center}
\caption{Estimates for the optimal $\alpha$ and $m$ (with a budget $p=40$) obtained via analytic minimization of \eqref{eq:mfmc_Jq} for the main (columns 2-5) and total (columns 6-9) Sobol index estimators. All values are averaged over 100 runs and estimated using 100 preprocessing samples. }
\label{table_ishigami_S}
\end{table}

\subsection{Quintic model}
\begin{figure}[!htb]
	\centering
	\includegraphics[width=.9\textwidth]{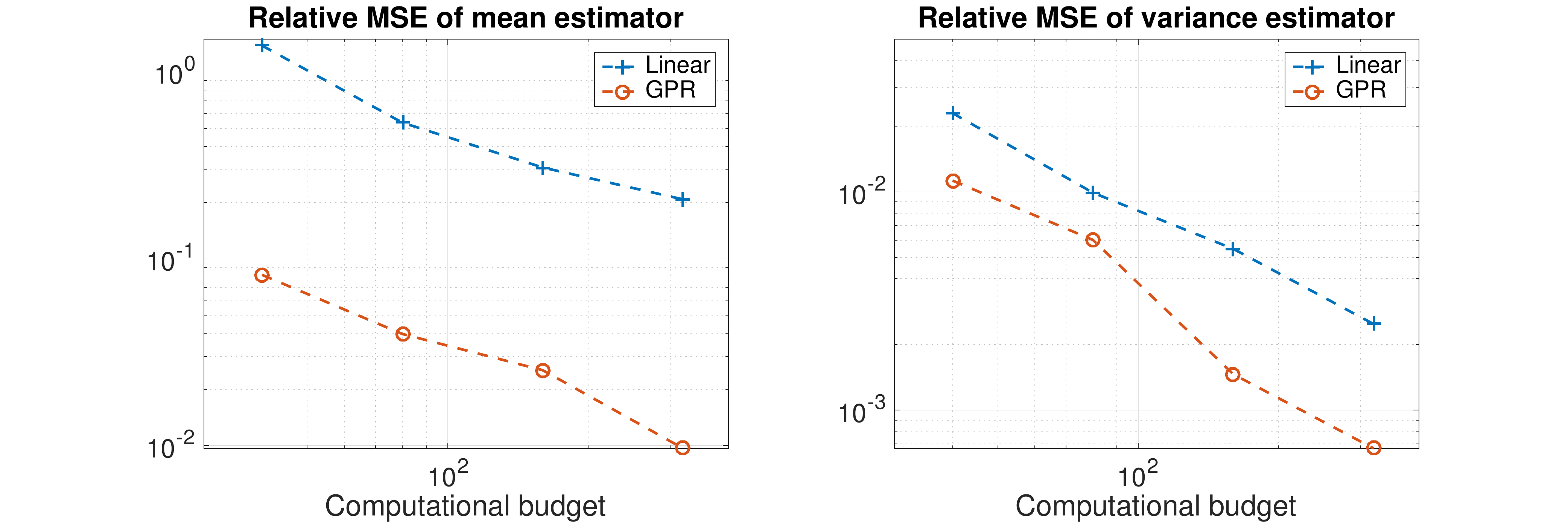}
  \caption{MSE of the proposed expectation and variance estimators for the Quintic model as a function of the number of computational budget $p$, computed with 100 replicates, for the linear and nonlinear estimators. Blue line: estimator using equation \eqref{eq:mfmc_q}. Red line: estimator using equation \eqref{eq:mfmc_g} and a GPR model is trained with 100 samples. }
	\label{quintic_MV}
\end{figure}

We here consider the following model
\begin{align*}
	f^{(1)} &=  \sin(z_1) +  \sin^2(z_2) +  \frac{1}{10} z_3^5, \quad z_i \sim \mathcal{U}(-\pi,\pi) \\
	f^{(2)} &=  \sin(z_1) +  \sin^2(z_2) +  2 z_3^3, \\
	f^{(3)} &=  \sin(z_1) +  \sin^2(z_2) + 20 z_3. 
\end{align*}
For this model, we estimated the exact expectation and variance using a standard Monte Carlo method. For the proposed multifidelity estimators of the expectation and variance, we use both approaches of linear correlation and statistical dependence, which are referred to in Figures \ref{quintic_MV} as linear and GPR, respectively. For both strategies, we also give the optimal budget allocations in Tables \ref{table_quintic_L} and \ref{table_quintic_NL}.

\begin{figure}[!htb]
	\centering
	\includegraphics[width=.48\textwidth]{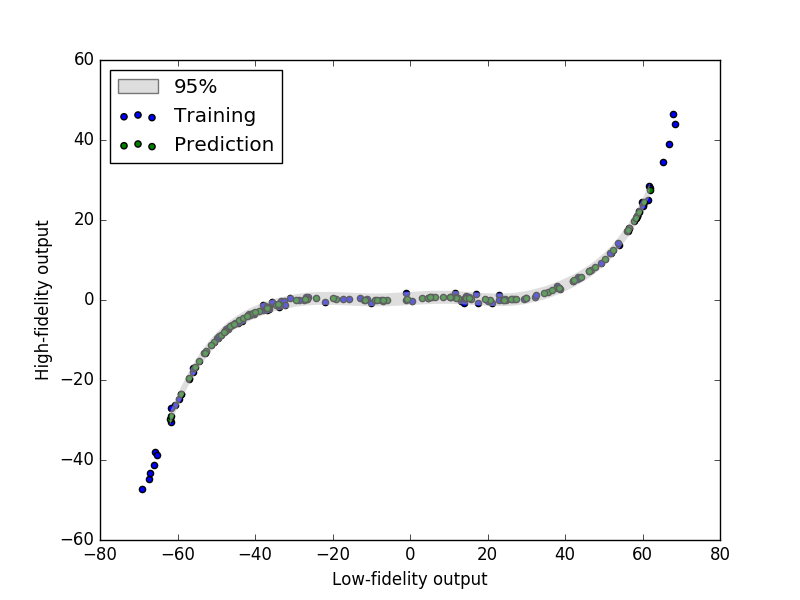}
	\includegraphics[width=.48\textwidth]{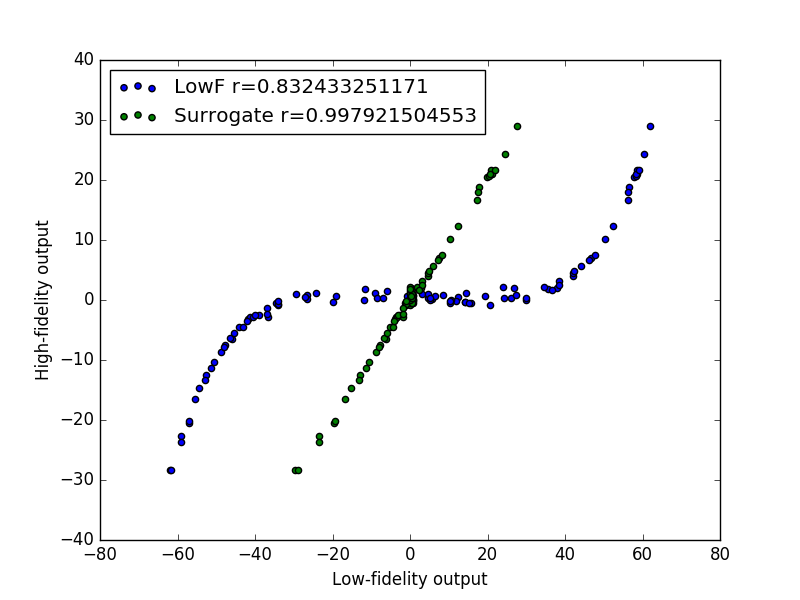}
  \caption{Left: Training (blue) and validation (green) used to train the GPR model mapping between the low-fidelity $f^{(2)}$ and high-fidelity $f^{(1)}$ models for the expectation of the Quintic problem, with confidence intervals of the prediction (gray). Right: Comparison of resulting linear correlations when using a linear (blue) or a nonlinear (green) function for mapping between the low-fidelity $f^{(2)}$ and high-fidelity $f^{(1)}$ models. }
	\label{quintic_correlation}
\end{figure}

\subsubsection{Expectation and variance estimation}
Tables \ref{table_quintic_L} and \ref{table_quintic_NL} summarize the values $\rho$ and the optimal $\alpha$ and $m$ for the case of, respectively, linear correlations and nonlinear statistical dependence. Both are averaged over 100 runs and the nonlinear model is estimated using 100 preprocessing samples. We see that in particular, the use of the nonlinear model significantly increases the correlation, to the point that the second model becomes redundant for estimating the expectation. Similarly, for the variance the number of model evaluations required by the high-fidelity reduces from 24 to 14 when the nonlinear dependence is exploited. Figure \ref{quintic_MV} shows the computed error of the expectation and variance estimators as a function of the budget, using 100 replicates of the multifidelity algorithm. From this chart, it is evident that taking advantage of the nonlinear statistical dependence, as done by the estimator \eqref{eq:mfmc_g}, provides a significant error reduction over the use of the more standard multifidelity estimator \eqref{eq:mfmc_q}, which is only based on linear correlations. 

\begin{table}[!htp]
\begin{center}
\begin{tabular}{lrrrrrr}
     	\toprule
     	\multicolumn{1}{l}{Models} &
     	\multicolumn{3}{c}{Expectation est.} &
     	\multicolumn{3}{c}{Variance est.} \\
     	\cmidrule(lr){2-4}
     	\cmidrule(lr){5-7}
	 & $m_k$ & $\alpha_k$ & $\rho_k$ & $m_k$ & $\alpha_k$ & $\rho_k$ \\
	 \midrule
	 $f^{(1)}$ & 25  & 1 & 1 & 24 & 1 & 1 \\
	 $f^{(2)}$ & 254 & 0.384 & 0.973 & 263 & 0.195 & 0.972  \\
	 $f^{(3)}$ & 2823 & 0.209 & 0.889 & 2521 & 0.119 & 0.784 \\
	 \bottomrule
\end{tabular}
\end{center}
\caption{Estimated optimal $\alpha$ and $m$ (with a budget $p=40$) of the Quintic model obtained via analytic minimization of \eqref{eq:mfmc_Jq} for the case of linear correlations. First three columns: values for the expectation estimator. Last three columns: values for the variance estimator. Both are averaged over 100 runs and estimated using 100 preprocessing samples. }
\label{table_quintic_L}
\end{table}

\begin{table}[!htp]
\begin{center}
\begin{tabular}{lrrrrrr}
     	\toprule
     	\multicolumn{1}{l}{Models} &
     	\multicolumn{3}{c}{Expectation est.} &
     	\multicolumn{3}{c}{Variance est.} \\
     	\cmidrule(lr){2-4}
     	\cmidrule(lr){5-7}
	 & $m_k$ & $\alpha_k$ & $\rho_k$ & $m_k$ & $\alpha_k$ & $\rho_k$ \\
	 \midrule
	 $f^{(1)}$ & 28  & 1 & 1 & 14 & 1 & 1 \\
	 $f^{(2)}$ & 0 & - & - & 417 & 0.971 & 0.996  \\
	 $f^{(3)}$ & 11768 &  0.99 & 0.997 & 5469 & 0.846 & 0.865 \\
	 \bottomrule
\end{tabular}
\end{center}
\caption{Estimated optimal $\alpha$ and $m$ (with a budget $p=40$) of the Quintic model obtained via analytic minimization of \eqref{eq:mfmc_Jg} for the case of nonlinear correlations. First three columns: values for the expectation estimator. Last three columns: values for the variance estimator. Both are averaged over 100 runs and estimated using 100 preprocessing samples. }
\label{table_quintic_NL}
\end{table}

\section{Application to cardiac electrophysiology}

We apply here the multifidelity methodology to a problem from
cardiac electrophysiology.  Briefly, the heart is an electrically
active organ, with an electric wave spreading throughout the tissue
thus to dictate the mechanical contraction.  The activation pattern,
defined as the first arrival time of the wave, is of clinical interest,
as an abnormal pattern leads to a mechanical dysfunction.
The propagation of the
traveling wave is anisotropic, being faster along the cardiac
fibers (usually 3- to 6-fold the transverse conduction velocity).
Fiber distribution is roughly known from histological studies~\cite{streeter1969}.
Patient-specific models incorporate the fiber distribution either
from diffusion-tensor MRI~\cite{Toussaint2013} or tailored to the anatomy
with a rule-based approach~\cite{potse2006,Bayer2012}. In both cases, epistemic
uncertainty is present.

In what follows, we consider the problem of quantifying the uncertainty
in the activation map, named $\AT(x)$, $x\in\Omega$ under the assumption
of randomly distributed cardiac fibers, $\vf(x)$.

The high-fidelity model is the monodomain equation~\cite{colli14}. This
is a reaction-diffusion equation, with nonlinear reaction, coupled
to a possibly large system of ODEs, encoding for the cellular
membrane model. The solution of the problem are the spatio-temporal
evolution of the transmembrane potential~$\Vm(x,t)$ and some auxiliary membrane
variables. The activation map is eventually obtained from the potential
as the time of first crossing of a given threshold:
\[
\AT_\text{h}(x) = \inf_t \bigl\{ t\colon \Vm(x,t) \ge V_\text{m,thres} \bigr\}.
\]
The monodomain equation reads as follows:
\[
\beta\Bigl( C_\text{m}\frac{\partial\Vm}{\partial t}
+ I_\text{ion}(\Vm,\mathbf{z}) - I_\text{stim}(x,t) \Bigr) = \nabla\cdot(\tGbar\nabla\Vm),
\]
where $C_\text{m}$ is the membrane capacitance, $\beta$ the surface-to-volume
ratio, $\mathbf{z}$ the set of auxiliary variables, and $\tGbar$ the
monodomain conductivity tensor, namely
\[
\tGbar(x) = \smt(x) \tI + \bigl(\sml(x)-\smt(x)\bigr) \vf(x)\otimes\vf(x),
\]
with $\sml$ and $\smt$ respectively the longitudinal and
transverse conductivity.

The low-fidelity model is the zeroth-order eikonal approximation of the
monodomain equation~\cite{pezzuto2017,colli14}.
The eikonal model is obtained by assuming that $\Vm(x,t)=U(t-\AT_\text{l}(x))$
for some action potential shape $U(\xi)$ and activation map $\AT_\text{l}$,
and then neglecting the diffusion term. The model reads:
\begin{equation}\label{eq:eiko}
\left\{ \begin{aligned}
& \frac{\theta}{\sqrt{\beta}}\sqrt{\tGbar\nabla\AT\cdot\nabla\AT} = 1,
& \mbox{in $\Omega \setminus \{ q_i \}_{i=1}^K$,} \\
& \AT(x_i) = \tau_i,
& i=1,\ldots,K.
\end{aligned}\right.
\end{equation}
The boundary conditions set the initial activation time to $\tau_i$
at some given locations $x_i\in\Omega$. The parameter~$\theta$ is a
scaling factor linked to the membrane model. We remark that for the
purposes of the analysis the function $U(\xi)$ is irrelevant.

In this study, we consider a patient-specific anatomy.
The electrophysiology parameters were already fitted to the patient data in
a previous study~\cite{potse2013}, with good
correlation in terms of activation map and surface ECG. The numerical grid
is composed by approximately 20 millions voxels at full resolution (0.2\,mm),
and 200\,000 voxels at coarse resolution (1\,mm).

For the monodomain simulation, we used Propag-5 software~\cite{Dickopf2014}.
Only the full resolution grid is considered, necessary to ensure
small numerical error~\cite{pezzuto2016,pathmanathan2014verification}. 
The membrane model is from \citet{tnnp04}. Activation is initiated with
a current stimulus in a 1\,mm$^3$ region for 2\,ms.
Regarding the eikonal model, we used our \textsc{GPGPU} implementation
described in~\cite{pezzuto2017}. We performed eikonal simulations at
both coarse and full resolution.



In summary, the full hierarchy of models for the multifidelity framework includes:
the discretized monodomain model, $f^{(1)}\colon\vf\mapsto\AT_\text{h}$, the discretized
eikonal model with high spatial resolution, $f^{(2)}\colon\vf\mapsto\AT_\text{l,1}$,
and the eikonal model at coarse spatial resolution,
$f^{(3)}\colon\vf\mapsto\AT_\text{l,2}$.

\subsection{Random fiber distribution and sampling}

The fiber distribution is modeled through a rule-based approach.
In details, the fiber orientation is locally determined by 3 angles, namely
$\alpha$, $\phi$ and $\gamma$. These angles are related to the Euler's angles,
but conveniently adapted so that $\alpha$ is consistent with the definition
given by~\citet{streeter1969}: it is the angle between the local fiber
direction and the valves' plane of the heart.
Transmurally, the angle $\alpha$ varies cubically from
$\pi/3$ at the endocardium to $-\pi/3$ at the epicardium~\cite{beyar1984}
\begin{equation}\label{eq:alpha0}
\alpha(r) = \frac{\pi}{3}(1-2r)^3,
\end{equation}
where $r$ measures the endocardium-to-epicardium distance, normalized to unity.
($r=0$ corresponds to the endocardium, $r=1$ to the epicardium.)

The random fiber field is obtained with the above rule, but assuming that
the angle $\alpha$ is a Gaussian random (scalar field), named $A(x)$, with
mean given by~\eqref{eq:alpha0} and correlation function as follows
\[
k(x,y) = \sigma\exp\Bigl(-\frac{d(x,y)^2}{2\ell^2}\Bigr),
\]
where $\sigma$ is the standard deviation, $\ell$ the correlation length,
and $d(x,y)$ the geodesic distance between $x$ and $y$.

The sampling strategy is provided in details in~\cite{Quaglino2018}. In
summary, the random angle $A(x)$ is approximated by a truncated Karhuhen-Lo\`eve
expansion:
\[
A(x) = A_0(x)
+ \sigma \sum^K_{i=1} \sqrt{\lambda_i} Z_i \varphi_i(x),
\]
with $Z_i$ independent standard normal random variables, and
$\lambda_i$ and $\varphi_i$ being eigenvalues and eigenvectors
of the Hilbert-Schmidt operator associated to the kernel $k(x,y)$.
The optimal truncation is obtained from a low-rank approximation
of the discretized correlation matrix
via a pivoted Cholesky decomposition~\cite{harbrecht2012}.
The distance function, required during the assembly of the
correlation matrix, is evaluated by solving an eikonal problem.

\begin{figure}[!htb]
	\centering
	\includegraphics[width=.95\textwidth]{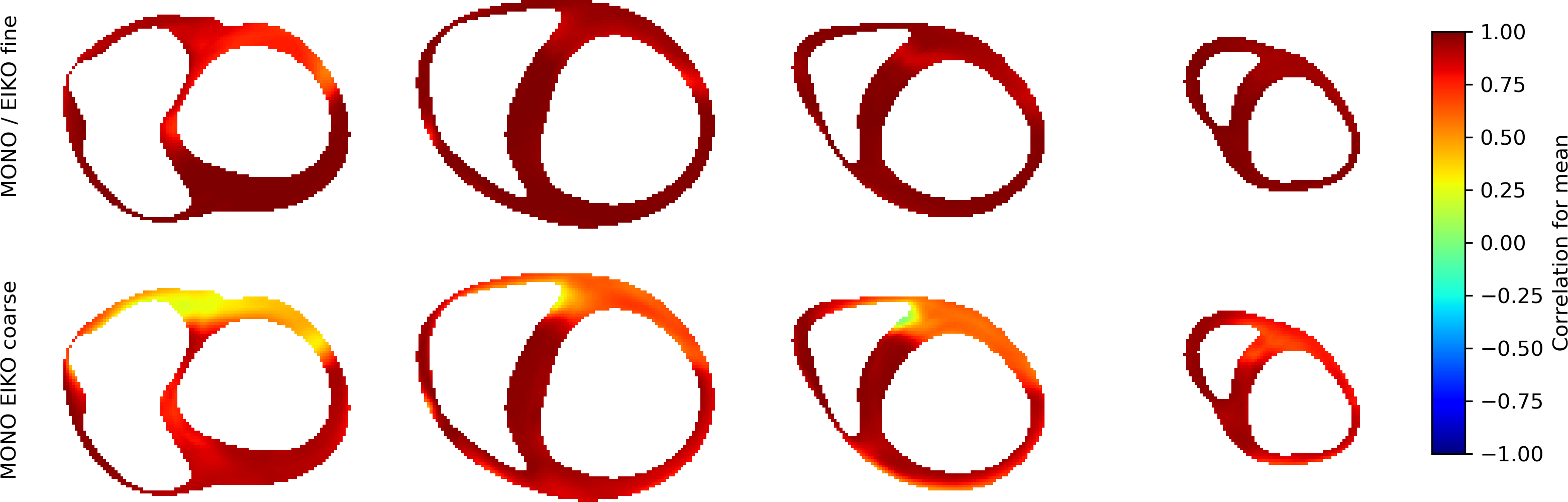}
  \caption{Correlations for the expectation estimator between monodomain and eikonal models
           at the same mesh resolution
           $h=0.02\,\si{cm}$ (top row) and at a coarser resolution $h=0.1\,\si{cm}$ (bottom row).
           Correlations are reported for slices of the heart from basal (left) to
           apical plane (right).}
	\label{patient_r}
\end{figure}

\begin{figure}[!htb]
	\centering
	\includegraphics[width=.95\textwidth]{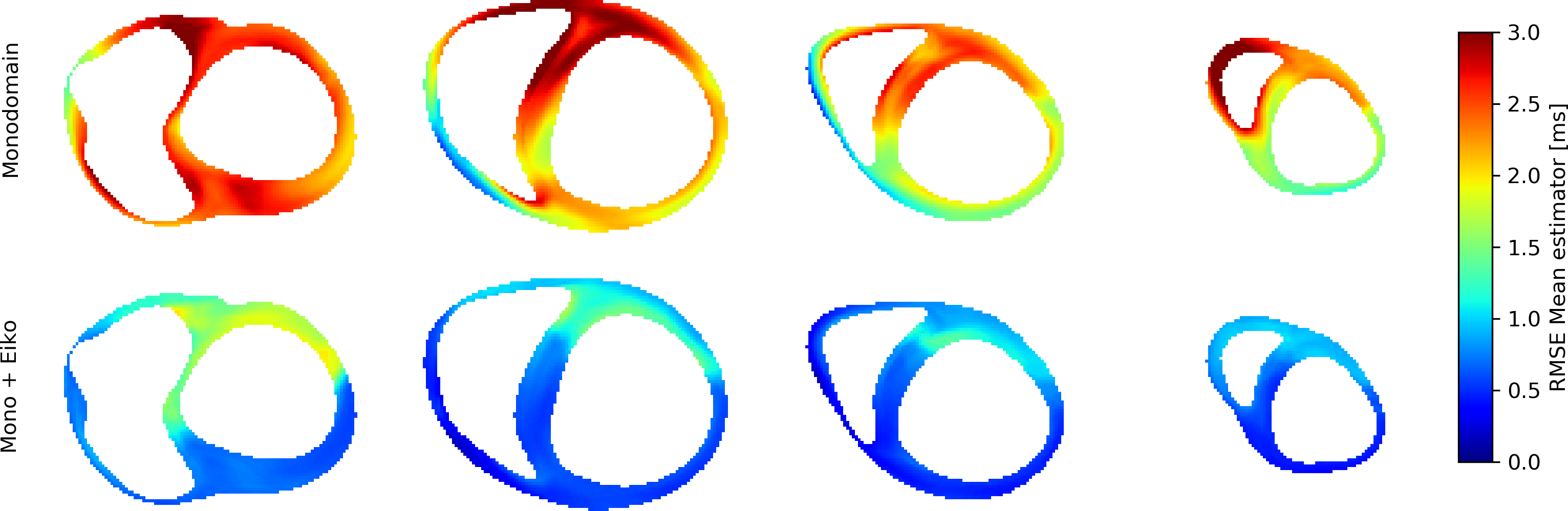}
  \caption{Estimated RMSE of the expectation estimator when only the monodomain model is used (top row),
            and its reduction when the eikonal model with fine grid is added (bottom row).
            Errors are reported for slices of the heart from basal plane (left) to
            apical plane (right).}
	\label{patient_e}
\end{figure}

\subsection{Expectation estimation}
Figure \ref{patient_r} shows the correlations between the expectation estimators for the bidomain and eikonal models at the same mesh resolution $h=0.02 \si{cm}$ (top) and at a coarser resolution $h=0.1 \si{cm}$ (bottom). Figure \ref{patient_e} shows the estimated error in the estimated expectation when only the bidomain model is used (top) and its reduction when the fine eikonal model is added (bottom). Figure \ref{patient_mean} (top) shows the computed expectation when all three models are used for the estimation. Table \ref{table_cardiac_V} summarizes the optimal number of model evaluations, estimated RMSEs $e_k$ with increasing number of models, and correlations $\bar{\rho}_k$ for a tolerance of $\epsilon=0.75~\si{ms}$. As a comparison, standard MC would require 92 samples of the high-fidelity model in order to achieve the same accuracy. The estimated node time for standard and multifidelity MC is, respectively, 277h and 36h.

\begin{table}[!htp]
\centering
\begin{tabular}{lrrrrrr}
     \toprule
     \multicolumn{1}{l}{Models} &
     \multicolumn{3}{c}{Expectation est.} &
     \multicolumn{3}{c}{Variance est.} \\
     \cmidrule(lr){2-4}
     \cmidrule(lr){5-7}
	 & $m_k$ & $e_k$ & $\bar{\rho}_k$ & $m_k$ & $e_k$ & $\bar{\rho}_k$ \\
	 \toprule
	 $f^{(1)}$ & 10  & 2.28~ms & 1 & 10 & 20.34~ms$^2$ & 1 \\
	 $f^{(2)}$ & 165 & 0.87~ms & 0.953 & 140 & 9.06~ms$^2$ &  0.905 \\
	 $f^{(3)}$ & 2490 & 0.76~ms & 0.750 & 1090 & 8.07~ms$^2$ & 0.567 \\
	 \bottomrule
\end{tabular}
\caption{Estimated optimal number of model evaluations $m$, RMSEs $e_k$ of the electrophysiology model obtained via analytic minimization of \eqref{eq:VarianceOptAlpha} with increasing number of models, and correlations $\bar{\rho}_k$. First two columns: values for the expectation estimator (tolerance of $\epsilon=0.75~\si{ms}$). Last two columns: values for the variance estimator (tolerance of $\epsilon=10.5~\si{ms^2}$). }
\label{table_cardiac_V}
\end{table}

\subsection{Variance estimation}

\begin{figure}[!htb]
	\centering
	\includegraphics[width=.95\textwidth]{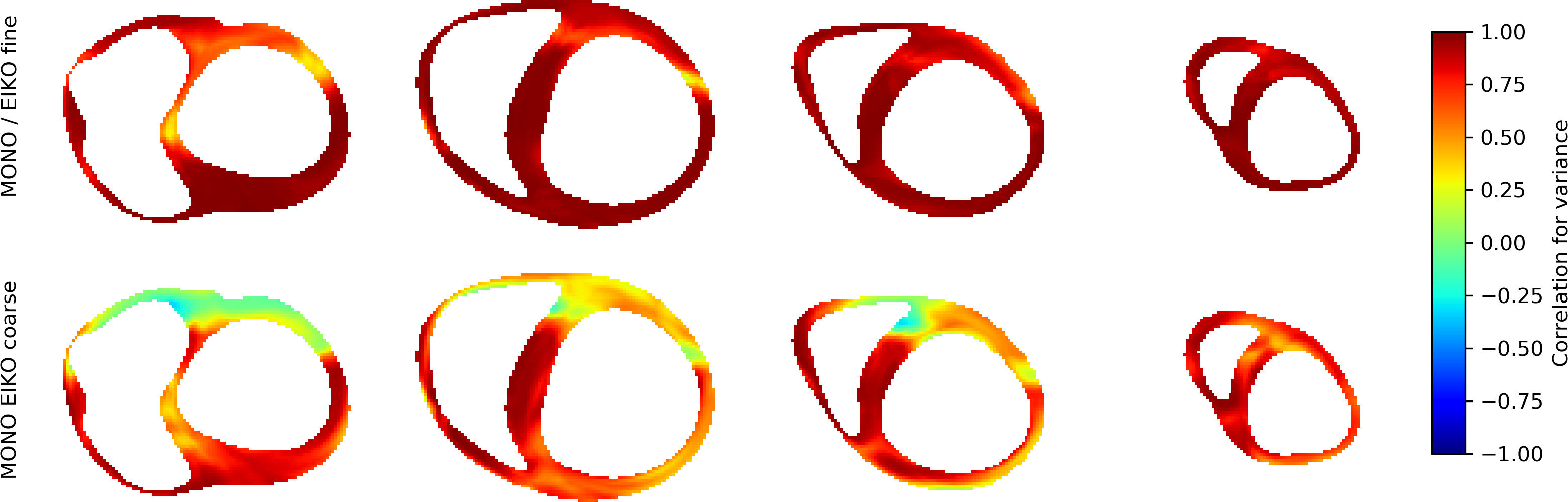}
  \caption{Correlations for the variance estimator between monodomain and eikonal models
           at the same mesh resolution
           $h=0.02\,\si{cm}$ (top row) and at a coarser resolution $h=0.1\,\si{cm}$ (bottom row).
           Correlations are reported for slices of the heart from basal (left) to
           apical plane (right).}
\label{patient_variance_r}
\end{figure}

\begin{figure}[!htb]
	\centering
	\includegraphics[width=.95\textwidth]{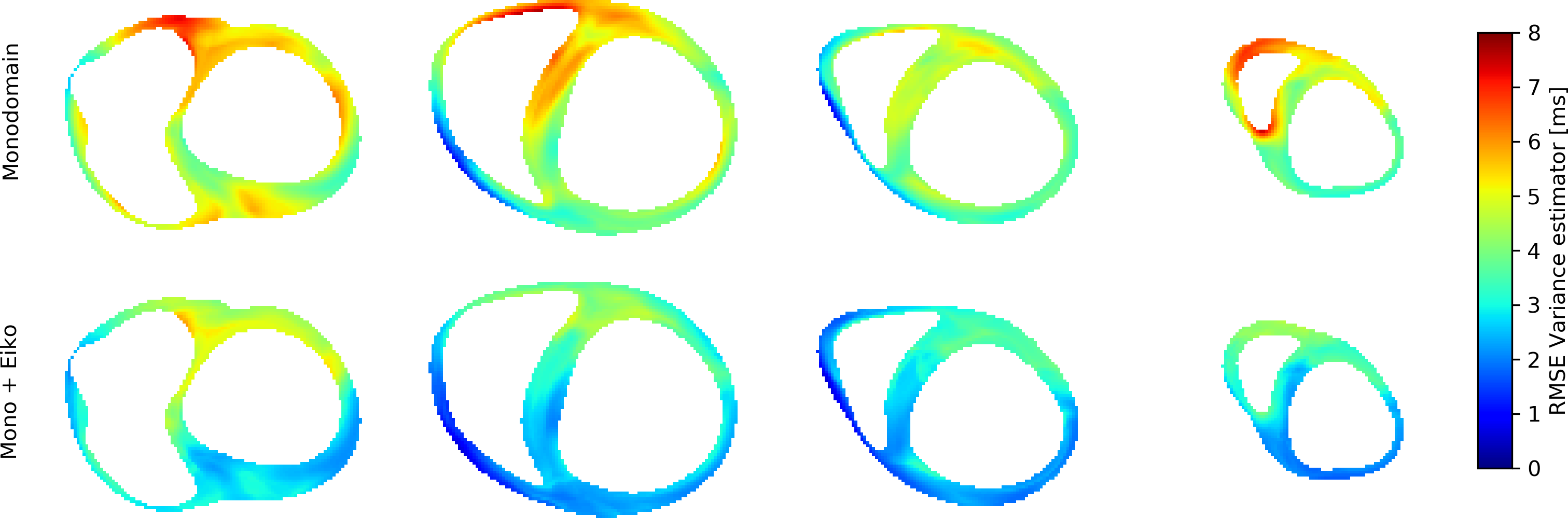}
  \caption{Estimated RMSE of the variance estimator when only the monodomal model is used (top row),
            and its reduction when the eikonal model with fine grid is added (bottom row).
            Errors are reported for slices of the heart from basal plane (left) to
            apical plane (right).}
	\label{patient_variance_e}
\end{figure}

\begin{figure}[!htb]
	\centering
	\includegraphics[width=.95\textwidth]{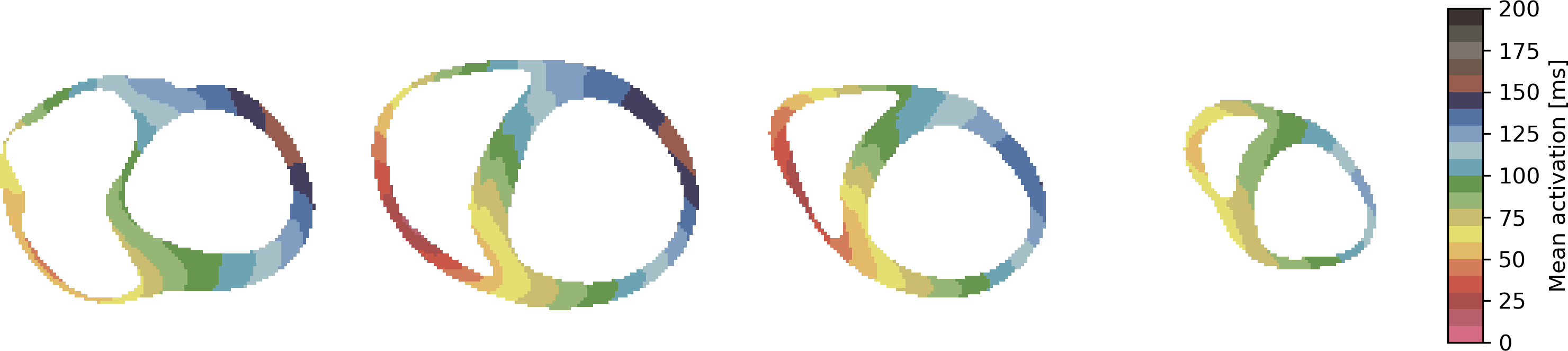}\\
	\includegraphics[width=.9325\textwidth]{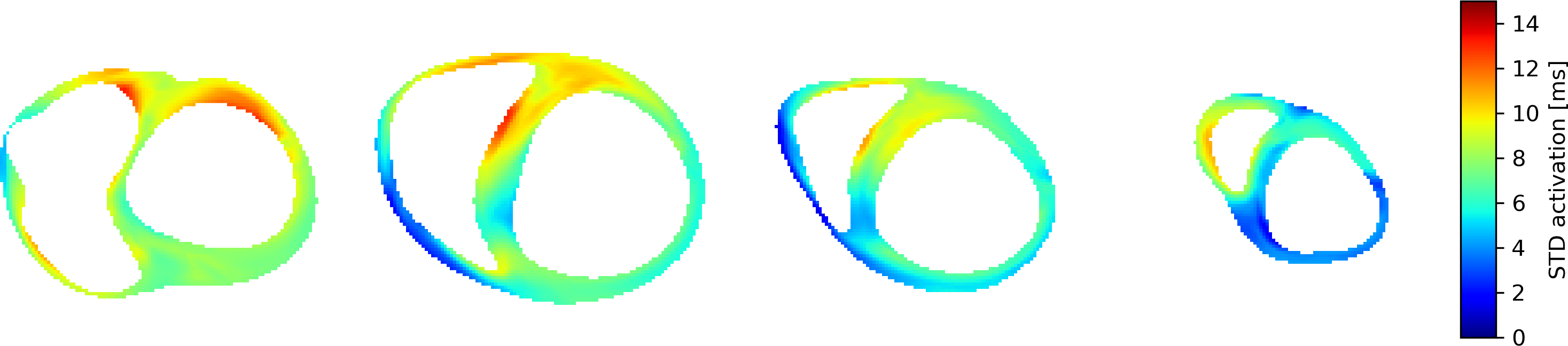}
  \caption{Estimated average of the activation map (top row) and estimated standard deviation (bottom row) when all when all three models are used.}
	\label{patient_mean}
\end{figure}

Figure \ref{patient_variance_r} shows the correlations between the variance estimators for the bidomain and eikonal models at the same mesh resolution $h=0.02 \si{cm}$ (top) and at a coarser resolution $h=0.1 \si{cm}$ (bottom). Figure \ref{patient_variance_e} shows the estimated error in the estimated standard deviation when only the bidomain model is used (top) and its reduction when the fine eikonal model is added (bottom). Figure \ref{patient_mean} (bottom) shows the computed standard deviation when all three models are used for the estimation. Table \ref{table_cardiac_V} summarizes the optimal number of model evaluations, estimated RMSEs $e_k$ with increasing number of models, and correlations $\bar{\rho}_k$ for a tolerance of $\epsilon=10.5~\si{ms^2}$. As a comparison, standard MC would require 48 samples of the high-fidelity model in order to achieve the same accuracy. The estimated node time for standard and multifidelity MC is, respectively, 145h and 35h.

%
%
%


\section{Application to nonlinear mechanics}

The last application is inspired by to the \emph{indentation problem}
of hyperelastic material, as described for instance in [REF]. Briefly, we
consider a $2\times2\times1\,\mbox{mm}^3$ cube clamped at the bottom,
while top face is indented of $0.1\,\mbox{mm}$ 
in a $1\,\mbox{mm}^2$ square region at its center.
We consider the finite deformation regime of a incompressible neo-Hookean
material with strain energy density function given as follows:
\[
\mathcal{W}(\mathbf{F}) = \frac{\mu}{2}(\mathbf{F}:\mathbf{F}-3) - \mu \ln(J).
\]
The deformation gradient tensor is defined as $\mathbf{F}=\mathbf{I}+\nabla\mathbf{u}$,
where $\mathbf{u}\colon\Omega\to\mathbb{R}^3$ is the (unknown) displacement.
The shear modulus is set to vary in space following a spatially-correlated
random field, evaluated via Karhuhen-Lo\`eve expantion:
\[
\mu(x) = \mu_0 + \sigma \sum_{i=1}^K \sqrt{\lambda_i} Z_i \phi_i(x),\qquad Z_i \sim \mathcal{N}(0,1).
\]
Functions $\phi_i(x)$ and coefficients $\lambda_i$ are, respectively,
eigenfunctions and eigenvalues of the Hilbert-Schmidt operator for
the squared-exponential correlation function, as in the previous section.
The correlation length is $\ell = \frac{1}{2\sqrt{2}}$, while mean and
standard deviation of $\mu$ are set to $80\,\mbox{MPa}$ and $20\,\mbox{MPa}$.
Samples with non-positive $\mu(x)$ are discarded.

The aim of this experiment is to estimated the pointwise average of the
von Mises stress, which is related to the second invariant of the Cauchy
stress:
\[
S = \frac{\mu}{J} \sqrt{\tfrac{3}{2} \dev\mathbf{B}:\dev\mathbf{B}},
\qquad \mathbf{B}=\mathbf{F}\mathbf{F}^{\mathrm{T}},
\]
given the uncertainty in the shear modulus.

The high-fidelity model is the Finite Element (FE) discretization of the problem
in mixed form with Taylor-Hood $\mathbb{P}^2-\mathbb{P}^1$ elements. We introduce
the pressure field $p$ to enforce the incompressibility constraint $J=1$.
The discrete problem is then solved with the Newton's method and a direct solver
for the linearised system.

The low-fidelity problem is based on the \emph{compressible} version of the
neo-Hookean material,
\[
\mathcal{W}_\mathrm{lf}(\mathbf{F}) = \frac{\mu}{2}(\mathbf{F}:\mathbf{F}-3) - \mu \ln(J)
+ \frac{\lambda}{2}(J-1)^2,
\]
with $\lambda = 10\mu_0 = 800\,\mbox{MPa}$. We discretized the problem with linear
FE for the displacement. Since the Poisson's ratio is about $0.45$, we do not expect
locking for this setup. In addition, the linear system within each Newton's iteration
is solved via algebraic multigrid.

\subsection{Expectation estimation of Von Mises stresses}
\begin{figure}[htb]
\centering
\includegraphics[width=0.32\textwidth,clip,trim=0 80 0 0]{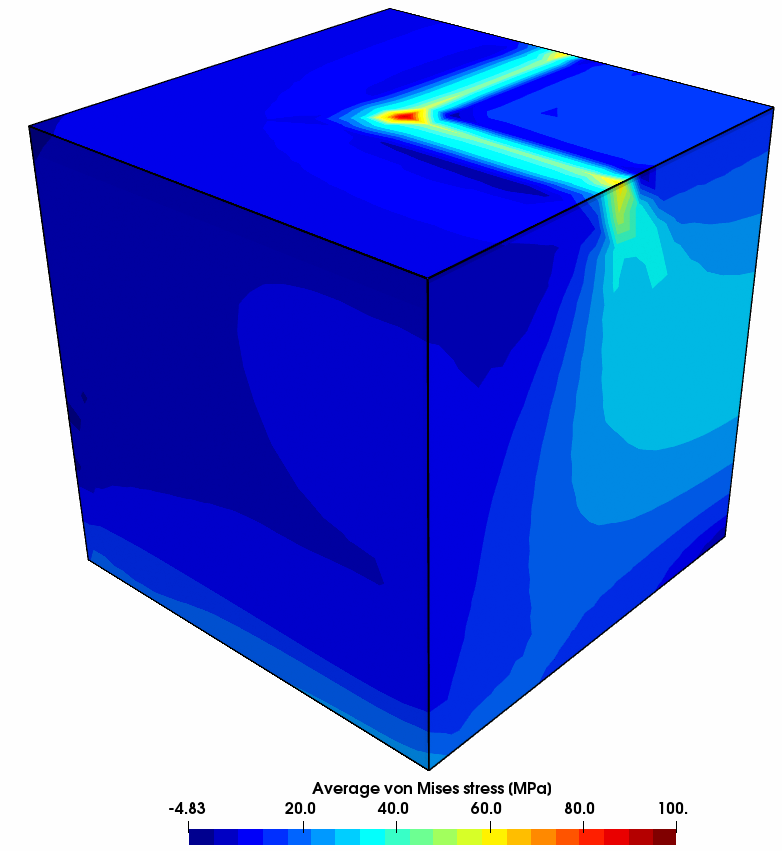}\hfill
\includegraphics[width=0.32\textwidth,clip,trim=0 80 0 0]{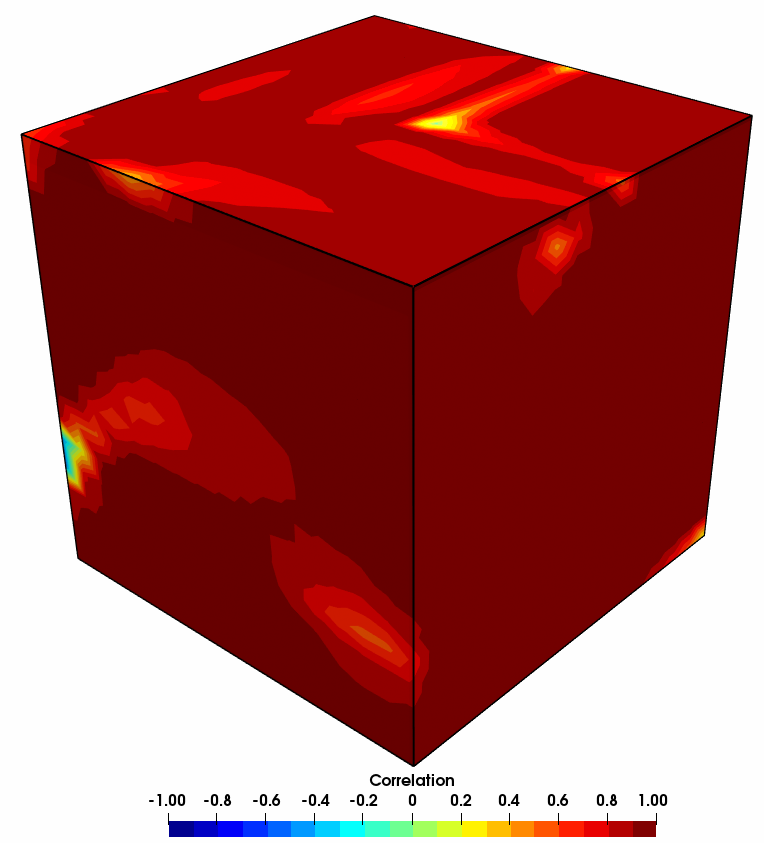}\hfill
\includegraphics[width=0.32\textwidth,clip,trim=0 80 0 0]{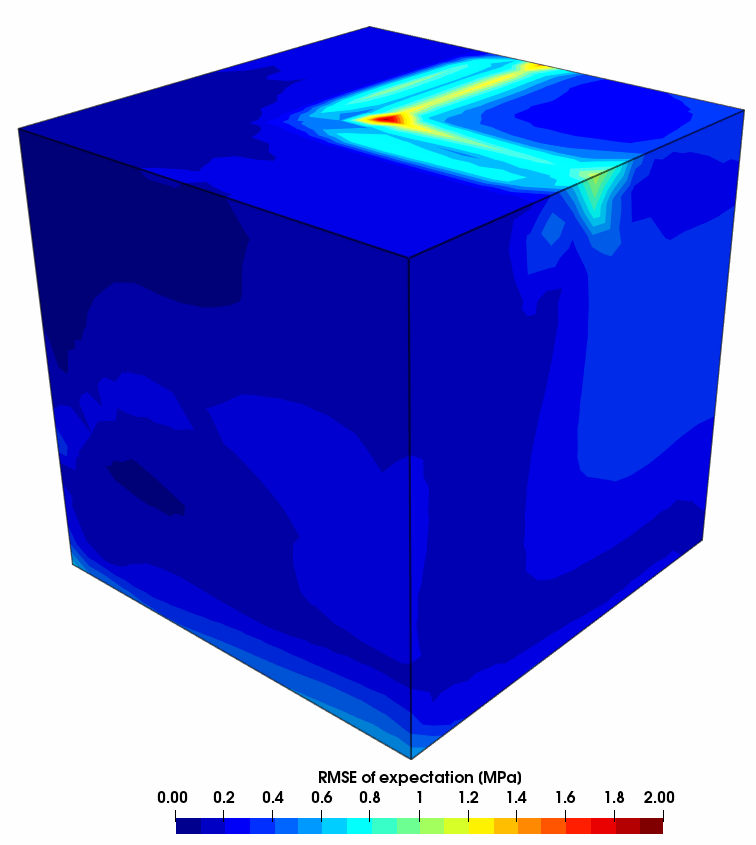}\\[1em]%
\includegraphics[width=0.32\textwidth,clip,trim=170 0 80 775]{mech_average}\hfill
\includegraphics[width=0.32\textwidth,clip,trim=170 0 80 775]{mech_corr}\hfill
\includegraphics[width=0.32\textwidth,clip,trim=170 5 80 770]{mech_error_mean}
\caption{MFMC estimate of average von Mises stress (left), pointwise correlation
coefficient (center), and RMSE of the estimator (right).}
\label{figure_mech}
\end{figure}

Figure \ref{figure_mech} (left) shows the computed expectation when both the high- and low-fidelity models are used for the estimation. Figure \ref{figure_mech} (center) shows the correlations between the expectation estimators for the two models. Figure \ref{figure_mech} (right) shows the estimated error in the estimated expectation when both models are used. We note that the high concentration of the error along the displaced boundary would have significantly increased the number of required model evaluations, had each point been treated as a separate scalar-valued QoI. Table \ref{table_mech} summarizes the optimal number of model evaluations, estimated RMSEs $e_k$ with increasing number of models, and correlations $\bar{\rho}_k$ for a budget of $B=2000~\si{s}$. As a comparison, standard MC would require 293 samples of the high-fidelity model in order to achieve the same accuracy. The estimated CPU time for standard and multifidelity MC is, respectively, $27000~\si{s}$ and $2000~\si{s}$.


\begin{table}[!htp]
\centering
\begin{tabular}{lrrrrrr}
     \toprule
     \multicolumn{1}{l}{Models} &
     \multicolumn{3}{c}{Expectation est.} &
     \multicolumn{3}{c}{Variance est.} \\
     \cmidrule(lr){2-4}
     \cmidrule(lr){5-7}
	 & $m_k$ & $e_k$ & $\bar{\rho}_k$ & $m_k$ & $e_k$ & $\bar{\rho}_k$ \\
	 \toprule
	 $f^{(1)}$ &  10 &    1.23~MPa     & 1 %
	           &  10 &   12.26~MPa$^2$ & 1 \\
	 $f^{(2)}$ & 550 &   0.223~MPa     & 0.992 %
	           & 395 &   2.85~MPa$^2$ & 0.985 \\
	 \bottomrule
\end{tabular}
\caption{Estimated optimal number of model evaluations $m$, RMSEs $e_k$ of the nonlinear mechanics model obtained via analytic minimization of \eqref{eq:VarianceOptAlpha} with increasing number of models, and correlations $\bar{\rho}_k$. First three columns: values for the expectation estimator
(budget $B=2000\,\mbox{s}$, $w_1=86\,\mbox{s}$, $w_2=1.85$\,\mbox{s}).
Last three columns: values for the variance estimator (same budget as previous set). }
\label{table_mech}
\end{table}

\subsection{Variance estimation of Von Mises stresses}
Table \ref{table_mech} (last three columns) summarizes the optimal number of model evaluations, estimated RMSEs $e_k$ with increasing number of models, and correlations $\bar{\rho}_k$ for a budget of $B=2000~\si{s}$. As a comparison, standard MC would require 184 samples of the high-fidelity model in order to achieve the same accuracy. The estimated CPU time for standard and multifidelity MC is, respectively, $16000~\si{s}$ and $2000~\si{s}$. The qualitative profiles of the estimated variance, the correlation, and the error, follow closely those of the expectation and are therefore not reported in the figures.


\section{Conclusions}
In this work, we extended and generalized previously proposed multifidelity
Monte Carlo estimators for the expectation \cite{peherstorfer2015optimal}, the variance, and sensitivity indices \cite{qian2017multifidelity}. We tested our  methodology on selected benchmarks, showing the good performance and usability of the proposed estimators. Moreover, we detailed their application to a real-world example from cardiac electrophysiology, which is unapproachable by standard Monte Carlo techniques.

Rather than considering the low- and high-fidelity models as opposing, which should be selected according to time constraints or accuracy considerations, the multifidelity framework combines them and exploits their best features. By doing so, a full UQ analysis of the complex problems, such as the monodomain equation in cardiac electrophysiology, can be performed in a total time comparable to that of using only the low-fidelity model, without sacrificing the complexity of the high-fidelity solution.

Future works will focus on reduced-physics models for electrophysiology, incorporating the repolarization and the ECG into the multifidelity hierarchy. Moreover, we intend to study multifidelity estimators for inverse problems.

\section*{Acknowledgements}
The authors acknowledge financial support by the Theo Rossi di Montelera Foundation, the Metis Foundation Sergio Mantegazza, the Fidinam Foundation, and the Horten Foundation to the Center for Computational Medicine in Cardiology. This work was also supported by grants from the Swiss National Supercomputing Centre (CSCS) under project ID s778.

\bibliography{mfmc}

\end{document}